\documentclass[10pt]{article}
\usepackage{amsmath,amsfonts,mathrsfs,amssymb}
\usepackage{indentfirst}

\begin{document}

\title{\textbf{Eigenvalue estimates of minimal hypersurfaces with finite index in Riemannian manifolds}}
\author{Zhongyang Sun$^{1,2}$\\
\small\footnotesize {1, $School$ $of$ $Mathematics$ $Science$, $Huaibei$ $Normal$ $University$, $Huaibei$, $Anhui$, 235000, $China$;}\\
\small\footnotesize {2, $Center$ $of$ $Mathematical$ $Sciences$, $Zhejiang$ $University$, $Hangzhou$, $Zhejiang$, 310027, $China$.}}
\date{}

\maketitle
\footnote[0]{2010 $Mathematics$ $Subject$ $Classification$. Primary 35P15; Secondary 53C40, 53C42.}

\footnote[0]{This work was supported by the National Natural Science Foundation of China (No. 11501500).}

\footnote[0]{E-mail: sunzhongyang12@163.com}

\vskip-.8cm {\normalsize $\mathbf{Abstract:}$ }{\small
The purpose of this paper is to study a
complete orientable minimal hypersurface with finite index
in an $(n+1)$-dimensional Riemannian manifold $N$.
We generalize Theorems  1.5-1.6 (\cite{Seo14}).
In 1976, Schoen and Yau proved the Liouville type theorem on stable minimal hypersurface, i.e., Theorem 1.7 (\cite{SchoenYau1976}).
Recently, Seo (\cite{Seo14}) generalized Theorem 1.7 (\cite{SchoenYau1976}). Finally, we generalize Theorems 1.7 (\cite{SchoenYau1976}) and 1.8 (\cite{Seo14}).
}

\vspace{6pt} {\normalsize $\mathbf{Keywords:}$} {\small
Stable minimal hypersurfaces; Riemannian manifolds ;
The first eigenvalue ; Hyperbolic space }

\section{Introduction and main results}

Cheng, Li, and Yau \cite{Cheng84} derived comparison theorems for the first eigenvalue of Dirichlet boundary problem on any compact domain in minimal submanifolds of the hyperbolic space by estimating the heat kernel of the compact
domain. Let $M^{n}$ be an $n$-dimensional smooth Riemannian manifold and $\Omega\in M^{n}$ a connected domain
with compact closure and nonempty boundary $\partial\Omega$.
The first eigenvalue of a Riemannian manifold $M^{n}$ is defined to be
$$
\lambda_{1}(\Omega)=\inf\left\{\frac{\int_{\Omega}|\nabla \phi|^{2}}{\int_{\Omega}\phi^{2}},~~\phi\in L_{1,0}^{2}(\Omega)\backslash\{0\}\right\},
$$
where the infimum is taken over all compactly supported Lipschitz functions on $M^{n}$ and $L_{1,0}^{2}(\Omega)$ is the completion of $C_{0}^{\infty}(\Omega)$
with respect to the norm
$$
\|\varphi\|_{\Omega}^{2}=\int_{\Omega}\varphi^{2}+\int_{\Omega}|\nabla\varphi|^{2}.
$$
If $\Omega_{1}\subset\Omega_{2}$ are bounded domains, then $0\leqslant\lambda_{1}(\Omega_{2})\leqslant\lambda_{1}(\Omega_{1})$.
Thus, the first eigenvalue of $M^{n}$ as the following limit
$$
\lambda_{1}(M^{n})=\lim_{r\to\infty}\lambda_{1}\big(B_{M^{n}}(p,r)\big)\geqslant0,
$$
where $B_{M^{n}}(p,r)$ is the geodesic ball of radius $r$ centered at $p$.

In 1970, McKean \cite{McKean70} studied the the bounded from below of the first eigenvalue of complete Riemannian manifold with sectional curvature bounded from above by negative constants and obtained the following famous result.

\vskip.2cm
{\noindent \bf Theorem 1.1 (\cite{McKean70})(McKean).}
\textit{Let $M^{n}$ be an $n$-dimensional complete simply connected Riemannian manifold with sectional curvature $K_{M^{n}}\leqslant-K_{1}^{2}<0$.
Then we have
$$
\lambda_{1}(M^{n})\geqslant\frac{(n-1)^{2}}{4}K_{1}^{2}.
$$}

In 2001,  Cheung-Leung \cite{CheungL01} gave a version of McKean's theorem for minimal submanifold of the hyperbolic space and obtained the following.

\vskip.2cm
{\noindent \bf Theorem 1.2 (\cite{CheungL01})(Cheung-Leung).}
\textit{Let $M^{n}(n\geqslant2)$ be a complete minimal submanifold in an $m$-dimensional hyperbolic space $\mathbb{H}^{m}(-1)$.
Then
$$
\lambda_{1}(M^{n})\geqslant\frac{(n-1)^{2}}{4}.
$$}

Recall that a minimal hypersurface is $stable$
if the second variation of the volume is nonnegative for any normal variation on a compact subset.
More precisely, a minimal hypersurface $M^{n}$ in an $(n+1)$-dimensional Riemannian manifold $\overline{M}^{n+1}$
is said to be called $stable$ if and only if for any $\phi\in C_{0}^{\infty}(M^{n})$
$$
\int_{M^{n}}\big[|\nabla \phi|^{2}-(\mathrm{\overline{Ric}}(\nu,\nu)+|A|^{2})\phi^{2}\big]d\upsilon\geqslant0,
$$
where $|A|^{2}$ is the squared length of the second fundamental form $A$, $\overline{\mathrm{Ric}}$ is the Ricci curvature of $\overline{M}^{n+1}$,
$\nu$ is the unit normal vector of $M^{n}$, and $d\upsilon$ is the volume form on ${M}^{n}$.

 In 2011, Seo \cite[Theorem 2.2]{Seo11} gave an upper bound for a stable minimal hypersurface with finite $L^{2}$-norm of the second fundamental form of $M^{n}$ and obtained the following.

\vskip.2cm
{\noindent \bf Theorem 1.3 (\cite{Seo11}).}
\textit{Let $M^{n}(n\geqslant2)$ be a complete stable minimal hypersurface in an $(n+1)$-dimensional hyperbolic space $\mathbb{H}^{n+1}(-1)$
with $\int_{M^{n}}|A|^{2}<\infty$.
Then
$$
\frac{(n-1)^{2}}{4}\leqslant\lambda_{1}(M^{n})\leqslant n^{2}.
$$}

Let $M^{n}$ be an $n$-dimensional complete manifold in an $(n+1)$-dimensional Riemannian manifold $N$.
We can choose a local vector field of orthonormal frames $e_{1},\cdots,e_{n+1}$ in $N$ such that the vectors $e_{1},\cdots,e_{n}$ are tangent to $M^{n}$ and the vector $e_{n+1}$ is normal to $M^{n}$.
In 2014, Seo generalized Theorem 1.3 (\cite{Seo11}) and obtain Theorem 1.4.

\vskip.2cm
{\noindent \bf Theorem  1.4 (\cite{Seo14}).}
\textit{Let $M^{n}$ be an $n$-dimensional complete non-totally geodesic stable minimal hypersurface
in an $(n+1)$-dimensional complete simply connected Riemannian manifold $N$ with sectional curvature $K_{N}$ satisfying that
$$
K_{1}\leqslant K_{N}\leqslant K_{2},
$$
where $K_{1}$, $K_{2}$ are constants and $K_{1}\leqslant K_{2}<0$.
Assume that, for $1-\sqrt{\frac{2}{n}}<p<1+\sqrt{\frac{2}{n}}$,
$$
\lim_{R\to\infty}R^{-2}\int_{B(R)}|A|^{2p}=0,
$$
where $B(R)$ is a geodesic ball of radius $R$ on $M^{n}$. If $|\nabla K_{N}|^{2}=\sum_{i,j,k,l,m}K_{ijkl;m}^{2}\leqslant K_{3}^{2}|A|^{2}$
for some constant $K_{3}\geqslant0$, where $K_{ijkl}$ is a curvature tensor of $N$ and $K_{ijkl;m}$ is the covariant derivative of $K_{ijkl}$.
Then we have
$$
-\frac{(n-1)^{2}}{4}K_{2}\leqslant \lambda_{1}(M^{n})\leqslant \frac{np^{2}\big(2K_{3}-n(K_{1}+K_{2})\big)}{2-n(p-1)^{2}}.
$$}

In this paper, first we generalize Theorem 1.4 and obtain Theorem A.

\vskip.2cm
{\noindent \bf Theorem  A.}
\textit{Let $M^{n}$ be an $n$-dimensional complete orientable non-totally geodesic minimal hypersurface with finite index
in an $(n+1)$-dimensional complete simply connected Riemannian manifold $N$ with sectional curvature $K_{N}$ satisfying that
$$
K_{1}\leqslant K_{N}\leqslant K_{2},
$$
where $K_{1}$, $K_{2}$ are constants and $K_{1}\leqslant K_{2}<0$.
Assume that, for $1-\sqrt{\frac{2}{n}}<p<1+\sqrt{\frac{2}{n}}$,
$$
\lim_{R\to\infty}R^{-2}\int_{B(R)}|A|^{2p}=0,
$$
where $B(R)$ is a geodesic ball of radius $R$ on $M^{n}$. If $|\nabla K_{N}|^{2}=\sum_{i,j,k,l,m}K_{ijkl;m}^{2}\leqslant K_{3}^{2}|A|^{2}$
for some constant $K_{3}\geqslant0$, where $K_{ijkl}$ is a curvature tensor of $N$ and $K_{ijkl;m}$ is the covariant derivative of $K_{ijkl}$.
Then we have
$$
-\frac{(n-1)^{2}}{4}K_{2}\leqslant \lambda_{1}(M^{n})\leqslant \frac{np^{2}\big(2K_{3}-n(K_{1}+K_{2})\big)}{2-n(p-1)^{2}}.
$$}

When we take $n=2$ and $p=1$ in Theorem A, we obtain the following.

\vskip.2cm
{\noindent \bf Corollary A.}
\textit{Let $M^{n}$ be an $n$-dimensional complete orientable non-totally geodesic minimal hypersurface with finite index
in an $(n+1)$-dimensional complete simply connected hyperbolic space $H^{n+1}$.
If $\int_{M^{n}}|A|^{2}<\infty$, then
$$
\frac{(n-1)^{2}}{4}\leqslant \lambda_{1}(M^{n})\leqslant n^{2}.
$$}

In 2012, Seo \cite{DungSeo12} proved that if $M^{n}$ is an $n$-dimensional complete noncompact stable minimal hypersurface in a Riemannian manifold with sectional
curvature bounded below by a nonpositive constant with $-K_{1}(2n-1)(n-1)<\lambda_{1}(M^{n})$, there is no nontrivial $L^{2}$ harmonic 1-form on $M^{n}$.
Next Seo \cite{Seo14} generalized the result and obtain Theorem 1.5.

\vskip.2cm
{\noindent \bf Theorem 1.5 (\cite{Seo14}).}
\textit{Let $M^{n}$ be an $n$-dimensional complete noncompact stable minimal hypersurface
in an $(n+1)$-dimensional  complete Riemannian manifold $N$ with sectional curvature
$K_{N}$ satisfying $K_{N}\geqslant K_{1}$, where $K_{1}\leqslant0$ is a constant.
Assume that, for $0<p<\frac{n}{n-1}+\sqrt{2n}$,
$$
\lambda_{1}(M^{n})>-\frac{2n(n-1)^{2}p^{2}K_{1}}{2n-\big[(n-1)p-n\big]^{2}}.
$$
Then there is no nontrivial $L^{2p}$ harmonic 1-form on $M^{n}$.}

We generalize Theorem 1.5 and get Theorem B.
\vskip.2cm
{\noindent \bf Theorem B.}
\textit{Let $M^{n}$ be an $n$-dimensional complete noncompact orientable minimal hypersurface with finite index
in an $(n+1)$-dimensional  complete Riemannian manifold $N$ with sectional curvature
$K_{N}$ satisfying $K_{N}\geqslant K_{1}$, where $K_{1}\leqslant0$ is a constant.
Assume that, for $0<p<\frac{n}{n-1}+\sqrt{2n}$,
$$
\lambda_{1}(M^{n})>-\frac{2n(n-1)^{2}p^{2}K_{1}}{2n-\big[(n-1)p-n\big]^{2}}.
$$
Then there is no nontrivial $L^{2p}$ harmonic 1-form on $M^{n}$.}

Recall that a function $f$ on a Riemannian manifold $M^{n}$ has $finite$ $L^{p}$ $energy$ if $|\nabla f|\in L^{p}(M^{n})$.
In 1976, Schoen and Yau proved the Liouville type theorem on stable minimal hypersurface as follow.

\vskip.2cm
{\noindent \bf Theorem 1.6 (\cite{SchoenYau1976}).}
\textit{Let $M^{n}$ be a complete noncompact stable minimal hypersurface in a Riemannian manifold N with nonnegative sectional curvature.
If $f$ is a harmonic function on $M^{n}$ with finite $L^{2}$ energy, $f$ is constant.}

Seo generalized Theorem 1.6 (\cite{SchoenYau1976}) and get the following.

\vskip.2cm
{\noindent \bf Theorem 1.7 (\cite{Seo14}).}
\textit{Let $M^{n}$ be an $n$-dimensional complete noncompact stable minimal hypersurface
in an $(n+1)$-dimensional Riemannian manifold $N$ with sectional curvature
$K_{N}$ satisfying $K_{N}\geqslant K_{1}$, where $K_{1}\leqslant0$ is a constant.
Assume that, for $0<p<\frac{n}{n-1}+\sqrt{2n}$,
$$
\lambda_{1}(M^{n})>-\frac{2n(n-1)^{2}p^{2}K_{1}}{2n-\big[(n-1)p-n\big]^{2}}.
$$
Then there is no nontrivial harmonic function on $M^{n}$ with finite $L^{p}$ energy.}

Finally, from Theorem B we obtain the following Corollary B, which generalize Theorems 1.6 (\cite{SchoenYau1976}) and 1.7 (\cite{Seo14}).

\vskip.2cm
{\noindent \bf Corollary B.}
\textit{Let $M^{n}$ be an $n$-dimensional complete noncompact orientable minimal hypersurface with finite index
in an $(n+1)$-dimensional Riemannian manifold $N$ with sectional curvature
$K_{N}$ satisfying $K_{N}\geqslant K_{1}$, where $K_{1}\leqslant0$ is a constant.
Assume that, for $0<p<\frac{n}{n-1}+\sqrt{2n}$,
$$
\lambda_{1}(M^{n})>-\frac{2n(n-1)^{2}p^{2}K_{1}}{2n-\big[(n-1)p-n\big]^{2}}.
$$
Then there is no nontrivial harmonic function on $M^{n}$ with finite $L^{p}$ energy.}

\vskip.2cm
{\noindent \bf Remark 1.}
We know that a complete stable minimal hypersurface in Riemannian
manifold with sectional curvature bounded below by a nonpositive constant has index 0. Hence
Theorems A and B can be regarded as generalizations of Theorems 1.4-1.5.
In particular, if $N$ is the $(n+1)$-dimensional hyperbolic space $H^{n+1}$, one sees
that $K_{1}=K_{2}=-1$, and hence $|\nabla K_{N}|^{2}=0$, that is, $K_{3}=0$.
Theorem A also generalizes Theorems 1.3 and 1.4.

Let $N$ be a $3$-dimensional hyperbolic space $H^{3}$ in Corollary A, we know that the finite index condition can be omitted, since the finiteness of the $L^{2}$ norm of the second fundamental form implies that $M^{2}$ has finite index, which was proved by
B¨¦rard et al. \cite{PCarmo97}.

\vskip.2cm
{\noindent \bf Theorem 1.8 (\cite{PCarmo97}).}
\textit{Let $M^{2}$ be a complete minimal hypersurface with $\int_{M^{2}}|A|^{2}<\infty$ in an $3$-dimensional complete hyperbolic space $\mathbb{H}^{3}(-1)$.
Then $M^{2}$ has finite index.}

From Corollary A and Theorem 1.8, we obtain

\vskip.2cm
{\noindent \bf Corollary C.}
\textit{Let $M^{2}$ be a complete stable  minimal hypersurface with $\int_{M^{2}}|A|^{2}<\infty$ in an $3$-dimensional complete hyperbolic space $\mathbb{H}^{3}(-1)$.
Then
$$
\frac{1}{4}\leqslant \lambda_{1}(M^{n})\leqslant 4.
$$}

\section{Proof of Theorem A}

\vskip.2cm
{\noindent \bf Proof of Theorem A.}
In 2003, Bessa-Montenegro \cite[Corollary 4.4]{Bessa03} obtained that the first eigenvalue $\lambda_{1}(M^{n})$ of a complete minimal hypersurface $M^{n}$ in $N$ is bounded below by $-\frac{(n-1)^{2}}{4}K_{2}>0$. In the rest of the proof, we only prove the upper bound of the first eigenvalue $\lambda_{1}(M^{n})$.

Since $M^{n}$ has finite index, there exists a compact subset $\Omega\subset M^{n}$ such that $M^{n}\backslash\Omega$ is stable, i.e., for any compactly supported Lipschitz function $\phi$ on $M^{n}\backslash\Omega$,
$$
\int_{M^{n}\backslash\Omega}\big[|\nabla\phi|^{2}-\big(|A|^{2}+\overline{\mathrm{Ric}}(e_{n+1})\big)\phi^{2}\big]d\upsilon\geqslant0,
\eqno(2.1)
$$
where $|A|^{2}$ denotes the squared length of the second fundamental form on $M^{n}\backslash\Omega$ and $d\upsilon$ denotes the volume form for the induced metric on $M^{n}\backslash\Omega$. For some geodesic ball $B(R_{0})\subset M^{n}$ centered at $x\in M^{n}$ of radius $R_{0}$ containing the compact set $\Omega$,
we know that $M^{n}\backslash B(R_{0})$ is still stable. Without loss of generality, we can assume that $\Omega=B(R_{0})$.

The assumption on sectional curvature of $N$ implies that
$$
nK_{1}\leqslant\overline{\mathrm{Ric}}(e_{n+1})=R_{(n+1)1(n+1)1}+\cdots+R_{(n+1)n(n+1)n}.
$$
Thus, the stability inequality (2.1) becomes
$$
\int_{M^{n}\backslash B(R_{0})}|\nabla\phi|^{2}-\big(|A|^{2}+nK_{1}\big)\phi^{2}\geqslant0.
\eqno(2.2)
$$
Choose a geodesic ball $B(R)\subset M^{n}$ centered at $x\in M^{n}$ of radius $R>R_{0}$ and take a cut-off function $0\leqslant f\leqslant 1$ on $M^{n}$ satisfying
$$
f=
\left\{
  \begin{array}{ll}
    0~~~~~~~~~\mathrm{on}~~~B(R_{0}),\\
    1~~~~~~~~~\mathrm{on}~~~B(2R+R_{0})\backslash B(R+R_{0}),\\
    0~~~~~~~~~\mathrm{on}~~~M^{n}\backslash B(3R+R_{0}),\\
  \end{array}
\right.
$$
and $|\nabla f|\leqslant \frac{1}{R}$ on $M^{n}$.

From the definition of $\lambda_{1}(M^{n})$ and the domain monotonicity of eigenvalue, it follows
$$
\lambda_{1}(M^{n})\leqslant\lambda_{1}\big(M^{n}\backslash B(R_{0})\big)\leqslant\frac{\int_{M^{n}\backslash B(R_{0})}|\nabla \phi|^{2}}{\int_{M^{n}\backslash B(R_{0})}\phi^{2}}
\eqno(2.3)
$$
for any $\phi\in W_{0}^{1,2}\big(M^{n}\backslash B(R_{0})\big)$, where $B(R_{0})\subset M^{n}$ is a geodesic ball  centered at $x\in M^{n}$ of radius $R_{0}$.

When we take $\phi=f|A|^{p}$ in (2.3), for a positive number $p>0$ and any compactly supported Lipschitz function $f$ on $M^{n}\backslash B(R_{0})$, we get
$$
\begin{aligned}
\lambda_{1}&(M^{n})\int_{M^{n}\backslash B(R_{0})}f^{2}|A|^{2p}\\
&\leqslant \int_{M^{n}\backslash B(R_{0})}\big|\nabla(f|A|^{p})\big|^{2}\\
&=\int_{M^{n}\backslash B(R_{0})}f^{2}\big|\nabla|A|^{p}\big|^{2}+\int_{M^{n}\backslash B(R_{0})}|\nabla f|^{2}|A|^{2p}+2\int_{M^{n}\backslash B(R_{0})}f|A|^{p}<\nabla f,\nabla|A|^{p}>.
\end{aligned}
\eqno(2.4)
$$
By using Young's inequality, for any $\varepsilon>0$, we have
$$
2\int_{M^{n}\backslash B(R_{0})}f|A|^{p}<\nabla f,\nabla|A|^{p}>\leqslant \varepsilon\int_{M^{n}\backslash B(R_{0})}f^{2}\big|\nabla|A|^{p}\big|^{2}+\frac{1}{\varepsilon}\int_{M^{n}\backslash B(R_{0})}|\nabla f|^{2}|A|^{2p}.
\eqno(2.5)
$$
From (2.4) and (2.5), we have
$$
\lambda_{1}(M^{n})\int_{M^{n}\backslash B(R_{0})}f^{2}|A|^{2p}\leqslant (1+\varepsilon)\int_{M^{n}\backslash B(R_{0})}f^{2}\big|\nabla|A|^{p}\big|^{2}+\left(1+\frac{1}{\varepsilon}\right)\int_{M^{n}\backslash B(R_{0})}|\nabla f|^{2}|A|^{2p}.
\eqno(2.6)
$$
On the other hand, combining the equations (1.22) and (1.27) in \cite{Schoen75}, at all points where $|A|\neq0$, we obtain
$$
|A|\Delta|A|+2K_{3}|A|^{2}-n(2K_{2}-K_{1})|A|^{2}+|A|^{4}\geqslant \sum h_{ijk}^{2}-\big|\nabla|A|\big|^{2}.
\eqno(2.7)
$$
Since $K_{2}-K_{1}\geqslant0$, from (2.7) we get
$$
\begin{aligned}
|A|\Delta|A|+2K_{3}|A|^{2}-nK_{2}|A|^{2}+|A|^{4}&\geqslant\sum h_{ijk}^{2}-\big|\nabla|A|\big|^{2}\\
&=|\nabla A|^{2}-\big|\nabla|A|\big|^{2}.
\end{aligned}
\eqno(2.8)
$$
In 2005, Y. L. Xin \cite{Xin05} obtained the following Kato-type inequality:
$$
|\nabla A|^{2}-\big|\nabla|A|\big|^{2}\geqslant\frac{2}{n}\big|\nabla|A|\big|^{2}.
\eqno(2.9)
$$
From (2.8) and (2.9), we obtain
$$
|A|\Delta|A|+2K_{3}|A|^{2}-nK_{2}|A|^{2}+|A|^{4}\geqslant\frac{2}{n}\big|\nabla|A|\big|^{2}.
\eqno(2.10)
$$
For a positive number $p>0$, we obtain
$$
\begin{aligned}
|A|^{p}\Delta|A|^{p}&=|A|^{p}\mathrm{div}(\nabla|A|^{p})\\
&=|A|^{p}\mathrm{div}(p|A|^{p-1}\nabla|A|)\\
&=p(p-1)|A|^{2p-2}\big|\nabla|A|\big|^{2}+p|A|^{2p-1}\Delta|A|\\
&=\frac{p-1}{p}\big|\nabla|A|^{p}\big|^{2}+p|A|^{2p-2}|A|\Delta|A|.
\end{aligned}
\eqno(2.11)
$$
Combining (2.10) and (2.11), we obtain
$$
\begin{aligned}
|A|^{p}\Delta|A|^{p}&\geqslant\frac{p-1}{p}\big|\nabla|A|^{p}\big|^{2}+\frac{2p}{n}|A|^{2p-2}\big|\nabla|A|\big|^{2}-p|A|^{2p+2}-p(2K_{3}-nK_{2})|A|^{2p}\\
&=\frac{p-1}{p}\big|\nabla|A|^{p}\big|^{2}+\frac{2}{np}\big|\nabla|A|^{p}\big|^{2}-p|A|^{2p+2}-p(2K_{3}-nK_{2})|A|^{2p}
\end{aligned}
$$
or equivalently,
$$
|A|^{p}\Delta|A|^{p}+p(2K_{3}-nK_{2})|A|^{2p}+p|A|^{2p+2}\geqslant\left(1-\frac{n-2}{np}\right)\big|\nabla|A|^{p}\big|^{2}.
\eqno(2.12)
$$
Multiplying (2.12) both side by a compactly supported Lipschitz function  $f^{2}$ and integrating over $B(3R+R_{0})\backslash B(R_{0})$,
we obtain
$$
\begin{aligned}
\int_{M^{n}\backslash B(R_{0})}f^{2}|A|^{p}\Delta|A|^{p}&+p(2K_{3}-nK_{2})\int_{M^{n}\backslash B(R_{0})}f^{2}|A|^{2p}+p\int_{M^{n}\backslash B(R_{0})}f^{2}|A|^{2p+2}\\
&\geqslant\left(1-\frac{n-2}{np}\right)\int_{M^{n}\backslash B(R_{0})}f^{2}\big|\nabla|A|^{p}\big|^{2}.
\end{aligned}
\eqno(2.13)
$$
The divergence theorem gives
$$
\begin{aligned}
0&=\int_{M^{n}\backslash B(R_{0})}\mathrm{div}(f^{2}|A|^{p}\nabla|A|^{p})\\
&=\int_{M^{n}\backslash B(R_{0})}f^{2}|A|^{p}\Delta|A|^{p}+\int_{M^{n}\backslash B(R_{0})}f^{2}\big|\nabla|A|^{p}\big|^{2}+2\int_{M^{n}\backslash B(R_{0})}f|A|^{p}<\nabla f,\nabla|A|^{p}>.
\end{aligned}
\eqno(2.14)
$$
From (2.13) and (2.14), we obtain
$$
\begin{aligned}
p\int_{M^{n}\backslash B(R_{0})}&f^{2}|A|^{2p+2}+p(2K_{3}-nK_{2})\int_{M^{n}\backslash B(R_{0})}f^{2}|A|^{2p}-\int_{M^{n}\backslash B(R_{0})}f^{2}\big|\nabla|A|^{p}\big|^{2}\\
&-2\int_{M^{n}\backslash B(R_{0})}f|A|^{p}<\nabla|A|^{p},\nabla f>\geqslant \left(1-\frac{n-2}{np}\right)\int_{M^{n}\backslash B(R_{0})}f^{2}\big|\nabla|A|^{p}\big|^{2}.
\end{aligned}
\eqno(2.15)
$$
Replacing $\phi$ by $f|A|^{p}$ in (2.2), we get
$$
\int_{M^{n}\backslash B(R_{0})}\big|\nabla(f|A|^{p})\big|^{2}\geqslant \int_{M^{n}\backslash B(R_{0})}f^{2}|A|^{2p+2}+nK_{1}\int_{M^{n}\backslash B(R_{0})}f^{2}|A|^{2p}
$$
or equivalently,
$$
\begin{aligned}
\int_{M^{n}\backslash B(R_{0})}|\nabla f|^{2}|A|^{2p}&+\int_{M^{n}\backslash B(R_{0})}f^{2}\big|\nabla|A|^{p}\big|^{2}+2\int_{M^{n}\backslash B(R_{0})}f|A|^{p}<\nabla f,\nabla|A|^{p}>\\
&\geqslant \int_{M^{n}\backslash B(R_{0})}f^{2}|A|^{2p+2}+nK_{1}\int_{M^{n}\backslash B(R_{0})}f^{2}|A|^{2p}.
\end{aligned}
\eqno(2.16)
$$
Multiplying (2.16) by a positive number $p>0$ and combining the inequalities (2.15), we obtain
$$
\begin{aligned}
p&\big[2K_{3}-n(K_{1}+K_{2})\big]\int_{M^{n}\backslash B(R_{0})}f^{2}|A|^{2p}+p\int_{M^{n}\backslash B(R_{0})}|A|^{2p}|\nabla f|^{2}\\
&+(p-1)\int_{M^{n}\backslash B(R_{0})}f^{2}\big|\nabla|A|^{p}\big|^{2}+2(p-1)\int_{M^{n}\backslash B(R_{0})}f|A|^{p}<\nabla f,\nabla|A|^{p}>\\
&\geqslant\left(1-\frac{n-2}{np}\right)\int_{M^{n}\backslash B(R_{0})}f^{2}\big|\nabla|A|^{p}\big|^{2}.
\end{aligned}
\eqno(2.17)
$$
Plugging (2.5) and (2.6) into (2.17), we have
$$
\begin{aligned}
\frac{p}{\lambda_{1}(M^{n})}&(2K_{3}-nK_{1}-nK_{2})\left[(1+\varepsilon)\int_{M^{n}\backslash B(R_{0})}\right.f^{2}\big|\nabla|A|^{p}\big|^{2}\\
&+\left.\left(1+\frac{1}{\varepsilon}\right)\int_{M^{n}\backslash B(R_{0})}|\nabla f|^{2}|A|^{2p}\right]+p\int_{M^{n}\backslash B(R_{0})}|A|^{2p}|\nabla f|^{2}\\
&+(p-1)\int_{M^{n}\backslash B(R_{0})}f^{2}\big|\nabla|A|^{p}\big|^{2}\\
&+(p-1)\left[\varepsilon\int_{M^{n}\backslash B(R_{0})}f^{2}\big|\nabla|A|^{p}\big|^{2}+\frac{1}{\varepsilon}\int_{M^{n}\backslash B(R_{0})}|\nabla f|^{2}|A|^{2p}\right]\\
&\geqslant\left(1-\frac{n-2}{np}\right)\int_{M^{n}\backslash B(R_{0})}f^{2}\big|\nabla|A|^{p}\big|^{2}
\end{aligned}
$$
or equivalently,
$$
\begin{aligned}
&\left[1-\frac{n-2}{np}-\left(1+\varepsilon\right)\left(\frac{p}{\lambda_{1}(M^{n})}(2K_{3}-nK_{1}-nK_{2})+p-1\right)\right]\int_{M^{n}\backslash B(R_{0})}f^{2}\big|\nabla|A|^{p}\big|^{2}\\
&\leqslant \left[\left(1+\frac{1}{\varepsilon}\right)\left(\frac{p}{\lambda_{1}(M^{n})}(2K_{3}-nK_{1}-nK_{2})+p\right)-\frac{1}{\varepsilon}\right]\int_{M^{n}\backslash B(R_{0})}|A|^{2p}|\nabla f|^{2}.
\end{aligned}
\eqno(2.18)
$$
Argue by contradiction. Assume that
$$
\lambda_{1}(M^{n})>\frac{np^{2}\big(2K_{3}-n(K_{1}+K_{2})\big)}{2-n(p-1)^{2}}=\frac{p(2K_{3}-nK_{1}-nK_{2})}{1-\frac{n-2}{np}-(p-1)}.
$$
Note that $1-\sqrt{\frac{2}{n}}<p<1+\sqrt{\frac{2}{n}}$ is equivalent to $2-n(p-1)^{2}>0$.
Choose a sufficiently small $\varepsilon>0$ satisfying that
$$
1-\frac{n-2}{np}-\left(1+\varepsilon\right)\left[\frac{p}{\lambda_{1}(M^{n})}(2K_{3}-nK_{1}-nK_{2})+p-1\right]>0.
\eqno(2.19)
$$
Using the fact that $|\nabla f|\leqslant\frac{1}{R}$ by our choice of $f$ and growth condition on $\int_{M^{n}\backslash B(R_{0})}|A|^{2p}$,
from (2.18) and (2.19) we can conclude that, by letting $R\rightarrow\infty$, $\int_{M^{n}}\big|\nabla|A|^{p}\big|^{2}=0$ on $M^{n}\backslash B(R_{0})$,
which implies that $|A|$ is constant on $M^{n}\backslash B(R_{0})$. Since the volume of $M^{n}$ is infinite \cite{Wei03},
we know that $|A|\equiv0$ outside the compact subset $B(R_{0})$.
It follows from the maximum principle for minimal hypersurface $M^{n}$ in Riemannian manifold $N$ that $M^{n}$ must be totally geodesic, which is impossible by our assumption.
Therefore, we have
$$
\lambda_{1}(M^{n})\leqslant \frac{np^{2}\big(2K_{3}-n(K_{1}+K_{2})\big)}{2-n(p-1)^{2}}.
$$
\hfill$\Box$

\section{Proof of Theorem B.}

In this section, in order to prove Theorem B, we need the following lemmas.

\vskip.2cm
{\noindent \bf Lemma 1 (\cite{Leung92}).}
\textit{Let $M^{n}$ be an n-dimensional complete immersed minimal
hypersurface in a Riemannian manifold $N$. If all the sectional curvatures of $N$ are bounded below by a constant $K$, then
$$
\mathrm{Ric}\geqslant (n-1)K-\frac{n-1}{n}|A|^{2}.
$$}

\vskip.2cm
{\noindent \bf Lemma 2 (\cite{Wang01}).}
\textit{Let $\omega$ be a harmonic 1-form in an $n$-dimensional Riemannian manifold $M^{n}$. Then
$$
|\nabla\omega|^{2}-\big|\nabla|\omega|\big|^{2}\geqslant\frac{1}{n-1}\big|\nabla|\omega|\big|^{2}.
$$}

{\noindent \bf Lemma 3 (\cite{Seo14}).}
\textit{Let $M^{n}$ be an n-dimensional complete noncompact Riemannian
manifold with $\lambda_{1}(M^{n})>0$. Then $\mathrm{Vol}(M^{n})=\infty$.
}

In the case of submanifolds, Cheung and Leung \cite[Corollary 2.2]{CheungL98} proved that the
volume $\mathrm{Vol}\big(B_{p}(r)\big)$ of every complete noncompact submanifold $M^{n}$ in the Euclidean
or hyperbolic space is infinite under the assumption
that the mean curvature vector $H$ of $M^{n}$ is bounded in absolute value.

Next we will prove vanishing theorem for $L^{p}$ harmonic 1-form
on a complete noncompact minimal hypersurface with finite index.

\vskip.2cm
{\noindent \bf Proof of Theorem B.}
Since $M^{n}$ has finite index,
by taking the similar processing
as in the proof of Theorem A, we can arrive to the following two inequalities
$$
\int_{M^{n}\backslash B(R_{0})}|\nabla\phi|^{2}-\big(|A|^{2}+nK_{1}\big)\phi^{2}\geqslant0;
\eqno(3.1)
$$
$$
\lambda_{1}(M^{n})\leqslant\lambda_{1}(M^{n}\backslash B(R_{0}))\leqslant\frac{\int_{M^{n}\backslash B(R_{0})}|\nabla \phi|^{2}}{\int_{M^{n}\backslash B(R_{0})}\phi^{2}}
\eqno(3.2)
$$
for any $\phi\in W_{0}^{1,2}\big(M^{n}\backslash B(R_{0})\big)$.

Let $\omega$ be an $L^{2p}$ harmonic 1-form on $M^{n}$, i.e.,
$$
\Delta\omega=0~~~~~\mathrm{and }~~~~~~\int_{M^{n}}|\omega|^{2p}d\upsilon<\infty.
$$
In an abuse of notation, we will refer to a harmonic 1-form and its dual harmonic vector field
both by $\omega$. From Bochner formula, it follows
$$
\Delta|\omega|^{2}=2\big(|\nabla\omega|^{2}+\mathrm{Ric}(\omega,\omega)\big).
$$
On the other hand, one sees that
$$
\Delta|\omega|^{2}=2\big(|\omega|\Delta|\omega|+\big|\nabla|\omega|\big|^{2}\big).
$$
Thus, we obtain
$$
|\omega|\Delta|\omega|-\mathrm{Ric}(\omega,\omega)=|\nabla\omega|^{2}-\big|\nabla|\omega|\big|^{2}.
$$
Applying Lemma 1 and Lemma 2 yields
$$
|\omega|\Delta|\omega|+\frac{n-1}{n}|A|^{2}|\omega|^{2}-(n-1)K_{1}|\omega|^{2}\geqslant\frac{1}{n-1}\big|\nabla|\omega|\big|^{2}.
\eqno(3.3)
$$
For any positive number $p$, we have
$$
\begin{aligned}
|\omega|^{p}\Delta|\omega|^{p}&=|\omega|^{p}\mathrm{div}(\nabla|\omega|^{p})\\
&=|\omega|^{p}\mathrm{div}(p|\omega|^{p-1}\nabla|\omega|)\\
&=p(p-1)|\omega|^{2p-2}\big|\nabla|\omega|\big|^{2}+p|\omega|^{2p-1}\Delta|\omega|\\
&=\frac{p-1}{p}\big|\nabla|\omega|^{p}\big|^{2}+p|\omega|^{2p-2}|\omega|\Delta|\omega|.
\end{aligned}
\eqno(3.4)
$$
Plugging inequality (3.3) into (3.4), we have
$$
\begin{aligned}
|\omega|^{p}\Delta|\omega|^{p}+p(n-1)\left(\frac{|A|^{2}}{n}-K_{1}\right)|\omega|^{2p}\geqslant\left[1-\frac{1}{p}+\frac{1}{p(n-1)}\right]\big|\nabla|\omega|^{p}\big|^{2}.
\end{aligned}
\eqno(3.9)
$$
Multiplying both side by a Lipschitz function $f^{2}$ with compact support in $M^{n}\backslash B(R_{0})$ and
integrating (3.5) over $M^{n}\backslash B(R_{0})$, we get
$$
\begin{aligned}
\int_{M^{n}\backslash B(R_{0})}f^{2}|\omega|^{p}\Delta|\omega|^{p}+\frac{p(n-1)}{n}&\int_{M^{n}\backslash B(R_{0})}f^{2}|A|^{2}|\omega|^{2p}-p(n-1)K_{1}\int_{M^{n}\backslash B(R_{0})}f^{2}|\omega|^{2p}\\
&\geqslant\left[1-\frac{1}{p}+\frac{1}{p(n-1)}\right]\int_{M^{n}\backslash B(R_{0})}f^{2}\big|\nabla|\omega|^{p}\big|^{2}
\end{aligned}
\eqno(3.6)
$$
The divergence theorem gives
$$
\int_{M^{n}\backslash B(R_{0})}f^{2}|\omega|^{p}\Delta|\omega|^{p}=-\int_{M^{n}\backslash B(R_{0})}f^{2}\big|\nabla|\omega|^{p}\big|^{2}-2\int_{M^{n}\backslash B(R_{0})}f|\omega|^{p}<\nabla f,\nabla|\omega|^{p}>.
\eqno(3.7)
$$
From (3.6) and (3.7), we have
$$
\begin{aligned}
\left[1-\frac{1}{p}+\frac{1}{p(n-1)}\right]&\int_{M^{n}\backslash B(R_{0})}f^{2}\big|\nabla|\omega|^{p}\big|^{2}\\
\leqslant &\frac{p(n-1)}{n}\int_{M^{n}\backslash B(R_{0})}f^{2}|A|^{2}|\omega|^{2p}-p(n-1)K_{1}\int_{M^{n}\backslash B(R_{0})}f^{2}|\omega|^{2p}\\
&-\int_{M^{n}\backslash B(R_{0})}f^{2}\big|\nabla|\omega|^{p}\big|^{2}-2\int_{M^{n}\backslash B(R_{0})}f|\omega|^{p}<\nabla f,\nabla|\omega|^{p}>.
\end{aligned}
\eqno(3.8)
$$
Replacing $\phi$ by $f|\omega|^{p}$ in (3.1), we have
$$
\int_{M^{n}\backslash B(R_{0})}\big|\nabla(f|\omega|^{p})\big|^{2}\geqslant \int_{M^{n}\backslash B(R_{0})}f^{2}|\omega|^{2p+2}+nK_{1}\int_{M^{n}\backslash B(R_{0})}f^{2}|\omega|^{2p}
$$
or equivalently,
$$
\begin{aligned}
\int_{M^{n}\backslash B(R_{0})}f^{2}\big|\nabla|\omega|^{p}\big|^{2}+&\int_{M^{n}\backslash B(R_{0})}|\nabla f|^{2}|\omega|^{2p}+2\int_{M^{n}\backslash B(R_{0})}f|\omega|^{p}<\nabla f, \nabla|\omega|^{p}>\\
&\geqslant\int_{M^{n}\backslash B(R_{0})}f^{2}|A^{2}|\omega|^{2p}+nK_{1}\int_{M^{n}\backslash B(R_{0})}f^{2}|\omega|^{2p}.
\end{aligned}
\eqno(3.9)
$$
Combining the inequalities (3.8) and (3.9) gives
$$
\begin{aligned}
\left[1-\frac{1}{p}+\frac{1}{p(n-1)}\right]&\int_{M^{n}\backslash B(R_{0})}f^{2}\big|\nabla|\omega\big|^{p}|^{2}\leqslant \frac{p(n-1)}{n}
\left\{\int_{M^{n}\backslash B(R_{0})}\right.f^{2}\big|\nabla|\omega\big|^{p}|^{2}\\
&+\int_{M^{n}\backslash B(R_{0})}|\nabla f|^{2}|\omega|^{2p}+2\int_{M^{n}\backslash B(R_{0})}f|\omega|^{p}<\nabla f,\nabla|\omega|^{p}>\\
&-nK_{1}\left.\int_{M^{n}\backslash B(R_{0})}f^{2}|\omega|^{2p}\right\}-p(n-1)K_{1}\int_{M^{n}\backslash B(R_{0})}f^{2}|\omega|^{2p}\\
&-\int_{M^{n}\backslash B(R_{0})}f^{2}\big|\nabla|\omega|^{p}\big|^{2}-2\int_{M^{n}\backslash B(R_{0})}f|\omega|^{p}<\nabla f,\nabla|\omega|^{p}>
\end{aligned}
$$
or equivalently,
$$
\begin{aligned}
\left[1-\frac{1}{p}+\frac{1}{p(n-1)}\right]&\int_{M^{n}\backslash B(R_{0})}f^{2}\big|\nabla|\omega|^{p}\big|^{2}\\
\leqslant& \left[\frac{p(n-1)}{n}-1\right]\int_{M^{n}\backslash B(R_{0})}f^{2}\big|\nabla|\omega|^{p}\big|^{2}\\
&+\frac{p(n-1)}{n}\int_{M^{n}\backslash B(R_{0})}|\nabla f|^{2}|\omega|^{2p}-2p(n-1)K_{1}\int_{M^{n}\backslash B(R_{0})}f^{2}|\omega|^{2p}\\
&+2\left[\frac{p(n-1)}{n}-1\right]\int_{M^{n}\backslash B(R_{0})}f|\omega|^{p}<\nabla f,\nabla|\omega|^{p}>.
\end{aligned}
\eqno(3.10)
$$
Put $\phi=f|\omega|^{p}$ in the above inequality (3.2). Then
$$
\begin{aligned}
\lambda_{1}&(M^{n})\int_{M^{n}\backslash B(R_{0})}f^{2}|\omega|^{2p}\\
&\leqslant \int_{M^{n}\backslash B(R_{0})}\big|\nabla(f|\omega|^{p})\big|^{2}\\
&=\int_{M^{n}\backslash B(R_{0})}f^{2}\big|\nabla|\omega|^{p}\big|^{2}+\int_{M^{n}\backslash B(R_{0})}|\nabla f|^{2}|\omega|^{2p}+2\int_{M^{n}\backslash B(R_{0})}f|\omega|^{p}<\nabla f,\nabla|\omega|^{p}>.
\end{aligned}
\eqno(3.11)
$$
By using Young's inequality, for any $\varepsilon>0$, we see that
$$
2\int_{M^{n}\backslash B(R_{0})}f|\omega|^{p}<\nabla f,\nabla|\omega|^{p}>\leqslant \varepsilon\int_{M^{n}\backslash B(R_{0})}f^{2}\big|\nabla|\omega|^{p}\big|^{2}+\frac{1}{\varepsilon}\int_{M^{n}\backslash B(R_{0})}|\nabla f|^{2}|\omega|^{2p}.
\eqno(3.12)
$$
From (3.11) and (3.12), we have
$$
\lambda_{1}(M^{n})\int_{M\backslash B(R_{0})}f^{2}|\omega|^{2p}\leqslant (1+\varepsilon)\int_{M^{n}\backslash B(R_{0})}f^{2}\big|\nabla|\omega|^{p}\big|^{2}+\left(1+\frac{1}{\varepsilon}\right)\int_{M^{n}\backslash B(R_{0})}|\nabla f|^{2}|\omega|^{2p}.
\eqno(3.13)
$$
Plugging (3.12) and (3.13) into (3.10), we have
$$
\begin{aligned}
&\left[2-\frac{1}{p}+\frac{1}{p(n-1)}+\frac{2p(n-1)K_{1}}{\lambda_{1}(M^{n})}-\frac{p(n-1)}{n}-\varepsilon\left(\frac{p(n-1)}{n}-1-\frac{2p(n-1)K_{1}}{\lambda_{1}(M^{n})}\right)\right]\\
&\times\int_{M^{n}\backslash B(R_{0})}f^{2}\big|\nabla|\omega|^{p}\big|^{2}\\
&\leqslant\left[\frac{p(n-1)}{n}-\frac{2p(n-1)K_{1}}{\lambda_{1}(M^{n})}-\frac{1}{\varepsilon}\left(\frac{p(n-1)}{n}-1-\frac{2p(n-1)K_{1}}{\lambda_{1}(M^{n})}\right)\right]
\int_{M^{n}\backslash B(R_{0})}|\nabla f|^{2}|\omega|^{2p}.
\end{aligned}
\eqno(3.14)
$$
Since
$$
\lambda_{1}(M^{n})>\frac{-2p(n-1)K_{1}}{2-\frac{1}{p}+\frac{1}{p(n-1)}-\frac{p(n-1)}{n}}=\frac{-2n(n-1)^{2}p^{2}K_{1}}{2n-\big[(n-1)p-n\big]^{2}},
$$
by the hypothesis, one can choose a sufficiently small $\varepsilon>0$ satisfying that
$$
2-\frac{1}{p}+\frac{1}{p(n-1)}+\frac{2p(n-1)K_{1}}{\lambda_{1}(M^{n})}-\frac{p(n-1)}{n}-\varepsilon\left(\frac{p(n-1)}{n}-1-\frac{2p(n-1)K_{1}}{\lambda_{1}(M^{n})}\right)>0.
\eqno(3.15)
$$
Note that $\int_{M^{n}\backslash B(R_{0})}|\omega|^{2p}<\infty$, since $\omega$ is an $L^{2p}$ harmonic 1-form on $M^{n}$. Letting $R$ tend
to infinity, from (3.14) and (3.15), we obtain
$$
\int_{M^{n}\backslash B(R_{0})}\big|\nabla|\omega|^{p}\big|^{2}=0,
$$
which implies that $|\nabla\omega|\equiv0$. Hence $|\omega|\equiv$constant. From the assumption that $\lambda_{1}(M^{n})>0$ and
Lemma 3, we can conclude that $|\omega|\equiv0$ outside the compact subset $B(R_{0})$.
It follows from the maximum principle for minimal hypersurface $M^{n}$ in Riemannian manifold $N$ that $|\omega|\equiv0$ on $M^{n}$.
\hfill$\Box$

\section*{Acknowledgement}
The author would like to thank Professor Kefeng Liu and Professor Hongwei Xu for long-time encouragement.

\bibliographystyle{amsplain}

\end{document}